\documentclass[a4paper,12pt,twoside]{article}
\usepackage{amsmath,amssymb}
\usepackage[amsmath,thmmarks]{ntheorem}
\usepackage{graphicx}
\usepackage{new}
\theoremseparator{.}
\theorembodyfont{\upshape}
\newtheorem{defn}{Definition}[section]
\newtheorem{rem}[defn]{Remark}
\newtheorem{conv}[defn]{Convention}
\theorembodyfont{\slshape}
\newtheorem{thm}[defn]{Theorem}
\newtheorem{prop}[defn]{Proposition}
\newtheorem{lem}[defn]{Lemma}

\newtheorem{fact}[defn]{Fact}
\theoremstyle{nonumberplain}
\newtheorem{mainthm}{Main Theorem}
\newtheorem{claim}{Claim}
\theorembodyfont{\upshape}
\newtheorem{ack}{Acknowledgements}
\theoremsymbol{\rule{.4em}{.8em}}
\theoremheaderfont{\itshape}
\newtheorem{prf}{Proof}
\newcommand{\word}[1]{\textit{#1}}

\newcommand{\A}{\mathcal{A}}
\newcommand{\B}{\mathcal{B}}
\newcommand{\F}{\mathcal{F}}
\newcommand{\I}{\mathcal{I}}
\newcommand{\K}{\mathcal{K}}
\newcommand{\N}{\mathbb{N}}

\title{Residuality of Families of $\F_{\sigma}$ Sets}
\author{Shingo SAITO}
\begin{document}
\begin{affiliation}
 Department of Mathematics,
 University College London, Gower Street, London, WC1E 6BT, England, United Kingdom.

 Email address: \texttt{ssaito@math.ucl.ac.uk}

\end{affiliation}

\maketitle

\begin{abstract}
 We prove that two natural definitions of residuality of families of $\F_{\sigma}$ sets are equivalent.
 We make use of the Banach-Mazur game in the proof.
\end{abstract}

\section{Introduction}
Properties of a typical compact set in the Euclidean space are often discussed.
Here we say that a property $P$ is fulfilled by a typical compact set if
the set of all compact sets satisfying $P$ is residual in the space of all compact sets
endowed with the Hausdorff metric.
It is well-known that a typical compact set in the Euclidean space is Lebesgue null
(see \cite{Zam}, for example).
In this paper we consider what a typical $\F_{\sigma}$ set means,
namely we define residuality of families of $\F_{\sigma}$ sets.
To the best of the author's knowledge, there has been no definition of such residuality.

We shall work in a compact, dense-in-itself metric space $(X,\rho)$ throughout this article.
Without loss of generality, we may assume that $\rho(x,y)\leqq1$ for any $x,y\in X$.
An $\F_{\sigma}$ set means an $\F_{\sigma}$ subset of $X$,
and $\F_{\sigma}$ stands for the set of all $\F_{\sigma}$ sets.
Let $\K$ denote the set of all compact (or equivalently closed) subsets of $X$.
For $x\in X$ and $r>0$, the closed ball of centre $x$ and radius $r$ is denoted by $\bar{B}(x,r)$.
For $K\in\K$ and $r>0$, we put $K[r]=\bigcup_{x\in K}\bar{B}(x,r)$.
It is well-known that the Hausdorff metric $d$ makes $\K$ a compact metric space.
Here we define $d(K,\emptyset)=1$ for any nonempty set $K\in\K$.
Then for $K,L\in\K$ and $r\in(0,1)$,
we have $d(K,L)\leqq r$ if and only if $K\subset L[r]$ and $L\subset K[r]$,
even when either $K$ or $L$ is empty.

Giving $\F_{\sigma}$ a topology would suffice to define residuality of families of $\F_{\sigma}$ sets,
but no good topology on $\F_{\sigma}$ has been found so far.
Bearing in mind that each $\F_{\sigma}$ set is the union of a sequence in $\K$,
we look at the space of sequences in $\K$ instead.
Here we might worry whether we should restrict ourselves only to increasing sequences,
but our main theorem removes this concern.
Let us proceed to rigorous definitions.

\begin{conv}
 Every sequence begins with the term of subscript one
 and the set $\N$ of all positive integers does not contain zero.
\end{conv}

The set of all sequences of sets in $\K$ is denoted by $\K^{\N}$
and endowed with the product topology.
The closed subset $\K_{\nearrow}^{\N}$ of $\K^{\N}$ is defined
as the set of all increasing sequences:
\[
 \K_{\nearrow}^{\N}=\bigl\{\,(K_n)\in\K^{\N}\bigm|K_1\subset K_2\subset\cdots\,\bigr\}.
\]

\begin{defn}
 For a family $\F$ of $\F_{\sigma}$ sets, we put
 \[
  \K_{\F}^{\N}=\Biggl\{\,(K_n)\in\K^{\N}\Biggm|\bigcup_{n=1}^{\infty}K_n\in\F\,\Biggr\}.
 \]
 We say that $\F$ is \word{$\K^{\N}$-residual} if
 $\K_{\F}^{\N}$ is residual in $\K^{\N}$
 and that $\F$ is \word{$\K_{\nearrow}^{\N}$-residual} if
 $\K_{\F}^{\N}\cap K_{\nearrow}^{\N}$ is residual in $\K_{\nearrow}^{\N}$.
\end{defn}

Our main theorem asserts that these two notions of residuality agree with each other:

\begin{mainthm}
 A family of $\F_{\sigma}$ sets is $\K^{\N}$-residual
 if and only if it is $\K_{\nearrow}^{\N}$-residual.
\end{mainthm}

The equivalence seems to show the appropriateness of our definitions.
Moreover our definitions match the properties of a typical compact set mentioned at the beginning.
We prove a lemma before we state the precise relation.

\begin{lem}\label{lem:Kuratowski-Ulam}
 Let $Y$ be a second countable topological space and $Z$ a nonempty Baire space.
 Then a subset $A$ of $Y$ is residual if and only if $A\times Z$ is residual in $Y\times Z$.
\end{lem}

\begin{prf}
 It suffices to show that a subset $A$ of $Y$ is meagre if and only if $A\times Z$ is meagre in $Y\times Z$.

 Suppose that $A$ is meagre.
 Then there exist nowhere dense sets $A_1$, $A_2$, \ldots\
 such that $A=\bigcup_{n=1}^{\infty}A_n$.
 It is easy to see that $A_n\times Z$ is nowhere dense in $Y\times Z$ for every $n\in\N$.
 Thus $A\times Z=\bigcup_{n=1}^{\infty}(A_n\times Z)$ is meagre.

 Conversely suppose that $A\times Z$ is meagre.
 Then the Kuratowski-Ulam theorem shows that
 for every $z$ in a residual set in $Z$,
 the set $\{\,y\in Y\mid(y,z)\in A\times Z\,\}=A$ is meagre.
 Therefore $A$ is meagre since $Z$ is a nonempty Baire space.
\end{prf}

\begin{rem}
 We shall use this lemma for $Y=\K$ and $Z=\K^{\N}$ in the next proposition.
 In this situation, the `if' part can be replaced by the following lemma,
 which is Lemma~4.25 of \cite{P} by Phelps:
 \begin{quote}
  Let $M$ be a complete metric space, $Y$ a Hausdorff space and
 $f\colon M\longrightarrow Y$ a continuous open surjective mapping.
  If $G$ is the intersection of countably many dense open subsets of $M$,
  then its image $f(G)$ is residual in $Y$.
 \end{quote}
 Indeed it suffices to substitute $\K\times\K^{\N}$ for $M$, $\K$ for $Y$, and the first projection for $f$.
 In order to prove this lemma, Phelps used the Banach-Mazur game,
 which we shall look at from the next section onwards.
\end{rem}

\begin{prop}
 Let $\I$ be a $\sigma$-ideal on $X$.
 Then $\I\cap\K$ is residual in $\K$ if and only if $\I\cap\F_{\sigma}$ is $\K^{\N}$-residual.
\end{prop}

\begin{prf}
 Since
 \begin{align*}
  \Biggl\{\,(K_n)\in\K^{\N}\Biggm|\bigcup_{n=1}^{\infty}K_n\in\I\,\Biggr\}
  &=\bigl\{\,(K_n)\in\K^{\N}\bigm|\text{$K_n\in\I$ for every $n\in\N$}\,\bigr\}\\
  &=\bigcap_{n=1}^{\infty}
    \bigl(\underbrace{\K\times\cdots\times\K}_{\text{$n-1$ times}}
          \times(\I\cap\K)\times\K\times\K\times\cdots\bigr),
 \end{align*}
 we see that $\I\cap\F_{\sigma}$ is $\K^{\N}$-residual if and only if
 $(\I\cap\K)\times\K\times\K\times\cdots$ is residual in $\K^{\N}$.
 Lemma~\ref{lem:Kuratowski-Ulam} shows that
 this is equivalent to the condition that
 $\I\cap\K$ is residual in $\K$.
\end{prf}

This proposition shows, for example, that a typical $\F_{\sigma}$ subset
of the interval $[0,1]$ is null.

\begin{ack}
 The author expresses his deep gratitude to his supervisor Professor David Preiss
 for invaluable suggestions and a lot of encouragement.
 In addition he is grateful to Mr Tim Edwards and Mr Hiroki Kondo
 for their careful reading of the manuscript.
 He also acknowledges the financial support
 by a scholarship from Heiwa Nakajima Foundation and
 by the Overseas Research Students Award Scheme.
\end{ack}

\section{Banach-Mazur games}
It is known that we can grasp residuality in terms of the Banach-Mazur game.

\begin{defn}\label{defn:BMgame}
 Let $Y$ be a topological space, $S$ a subset of $Y$, and $\A$ a family of subsets of $Y$.
 Suppose that every set in $\A$ has nonempty interior
 and that every nonempty open subset of $Y$ contains a set in $\A$.
 The \word{$(Y,S,\A)$-Banach-Mazur game} is described as follows.
 Two players, called Player I and Player II, alternately choose a set in $\A$
 with the restriction that
 they must choose a subset of the set chosen in the previous turn.
 Player~II will win if the intersection of all the sets chosen by the players is contained in $S$;
 otherwise Player~I will win.
\end{defn}

\begin{rem}
 The assumptions on $\A$ ensure that the players can continue to take sets.
\end{rem}

\begin{fact}\label{fact:BMgame}
 The $(Y,S,\A)$-Banach-Mazur game has a winning strategy for Player~II
 if and only if $S$ is residual in $Y$.
\end{fact}

For the proof of this fact, we refer the reader to Theorem~1 in \cite{Ox}.

\bigbreak

In order to prove our main theorem, we look at the following Banach-Mazur games:

\begin{defn}\label{defn:K^N-BMgame}
 Let $\F$ be a family of $\F_{\sigma}$ sets.

 Let $\B$ denote the family of all sets of the form
 \[
  \bar{B}\bigl((K_n),a,r\bigr)=
  \bigl\{\,(A_n)\in\K^{\N}\bigm|\text{$d(A_n,K_n)\leqq r$ for $n=1,\ldots,a$}\,\bigr\},
 \]
 where $a$ is a positive integer,
 $(K_n)$ is a sequence in $\K^{\N}$ such that $K_1$, \ldots, $K_a$ are pairwise disjoint finite sets,
 and $r$ is a positive real number less than $1$ such that
 any two distinct points in $\bigcup_{j=1}^{a}K_j$ have distance at least $3r$.
 The $(\K^{\N},\K_{\F}^{\N},\B)$-Banach-Mazur game
 is called the \word{$(\K^{\N},\F)$-BM game} for ease of notation.

 Let $\B_{\nearrow}$ denote the family of all sets of the form
 \[
  \bar{B}_{\nearrow}\bigl((L_n),b,s\bigr)=
  \bigl\{\,(A_n)\in\K_{\nearrow}^{\N}\bigm|\text{$d(A_n,L_n)\leqq s$ for $n=1,\ldots,b$}\,\bigr\},
 \]
 where $b$ is a positive integer,
 $(L_n)$ is a sequence in $\K_{\nearrow}^{\N}$
 such that $L_1$, \ldots, $L_b$ are finite,
 and $s$ is a positive real number less than $1$ such that
 any two distinct points in $L_b$ have distance at least $3s$.
 The $(\K_{\nearrow}^{\N},\K_{\F}^{\N}\cap\K_{\nearrow}^{\N},\B_{\nearrow})$-Banach-Mazur game
 is called the \word{$(\K_{\nearrow}^{\N},\F)$-BM game}.
\end{defn}

\begin{rem}
 Notice that the families $\B$ and $\B_{\nearrow}$ satisfy the assumptions in Definition~\ref{defn:BMgame}
 since $X$ is dense-in-itself.
\end{rem}

\begin{conv}
 Whenever we write $\bar{B}\bigl((K_n),a,r\bigr)$ or $\bar{B}_{\nearrow}\bigl((L_n),b,s\bigr)$,
 we assume that
 $(K_n)$, $a$, $r$; $(L_n)$, $b$, $s$
 satisfy the conditions in Definition~\ref{defn:K^N-BMgame}.
\end{conv}

\begin{rem}
 A trivial observation shows that
 $\bar{B}\bigl((K_n),a,r\bigr)\subset\bar{B}\bigl((K_n'),a',r'\bigr)$ implies
 $a\geqq a'$ and $r\leqq r'$
 and that
 $\bar{B}_{\nearrow}\bigl((L_n),b,s\bigr)\subset\bar{B}_{\nearrow}\bigl((L_n'),b',s'\bigr)$ implies
 $b\geqq b'$ and $s\leqq s'$.
\end{rem}

Fact \ref{fact:BMgame} enables us to translate our main theorem into the following:

\begin{thm}\label{thm:main-BM}
 For a family $\F$ of $\F_{\sigma}$ sets,
 the $(\K^{\N},\F)$-BM game has a winning strategy for Player~II
 if and only if the $(\K_{\nearrow}^{\N},\F)$-BM game does.
\end{thm}

\section{Proof of our main theorem}
In this section we shall prove Theorem \ref{thm:main-BM},
which, as we have already mentioned, implies our main theorem.
Hereafter we fix a family $\F$ of $\F_{\sigma}$ sets
and call the Banach-Mazur games without referring to $\F$.

\subsection{Outline of the proof}
\begin{figure}
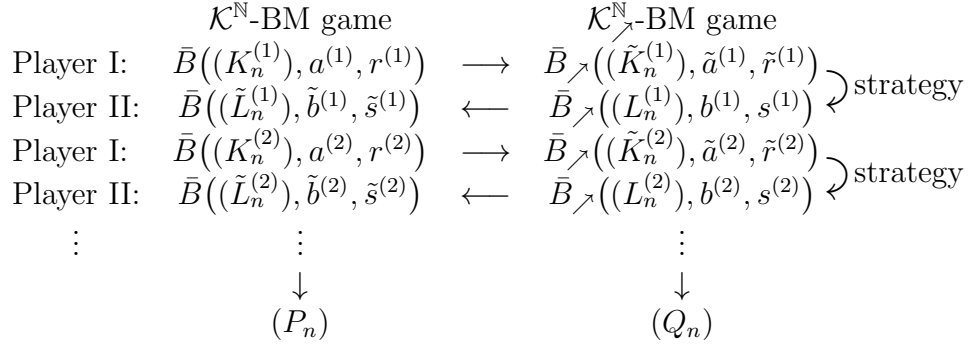

 \centering
 \begin{tabular}{lccc}
  &$\K^{\N}$-BM game&&$\K_{\nearrow}^{\N}$-BM game\\
  Player~I:&
  $\bar{B}\bigl((K_n^{(1)}),a^{(1)},r^{(1)}\bigr)$&$\longrightarrow$&%
  $\bar{B}_{\nearrow}\bigl((\tilde{K}_n^{(1)}),\tilde{a}^{(1)},\tilde{r}^{(1)}\bigr)$%
  \makebox[0pt][l]{\hspace*{2pt}\raisebox{2pt}[0pt][0pt]{\rotatebox{-90}{\Large$\curvearrowright$}}%
                   \hspace*{3pt}\raisebox{-8pt}[0pt][0pt]{strategy}}\\
  Player~II:&
  $\bar{B}\bigl((\tilde{L}_n^{(1)}),\tilde{b}^{(1)},\tilde{s}^{(1)}\bigr)$&$\longleftarrow$&%
  $\bar{B}_{\nearrow}\bigl((L_n^{(1)}),b^{(1)},s^{(1)}\bigr)$\\
  Player~I:&
  $\bar{B}\bigl((K_n^{(2)}),a^{(2)},r^{(2)}\bigr)$&$\longrightarrow$&%
  $\bar{B}_{\nearrow}\bigl((\tilde{K}_n^{(2)}),\tilde{a}^{(2)},\tilde{r}^{(2)}\bigr)$%
  \makebox[0pt][l]{\hspace*{2pt}\raisebox{2pt}[0pt][0pt]{\rotatebox{-90}{\Large$\curvearrowright$}}%
                   \hspace*{3pt}\raisebox{-8pt}[0pt][0pt]{strategy}}\\
  Player~II:&
  $\bar{B}\bigl((\tilde{L}_n^{(2)}),\tilde{b}^{(2)},\tilde{s}^{(2)}\bigr)$&$\longleftarrow$&%
  $\bar{B}_{\nearrow}\bigl((L_n^{(2)}),b^{(2)},s^{(2)}\bigr)$\\
  \multicolumn{1}{c}{\vdots}&\vdots&&\vdots\\
  &$\downarrow$&&$\downarrow$\\
  &$(P_n)$&&$(Q_n)$
 \end{tabular}
 \caption{Outline of the proof that $\K_{\nearrow}^{\N}$-residuality implies $\K^{\N}$-residuality}
 \label{fig:transfer_K^N-->K_nearrow^N}
\end{figure}

\begin{figure}
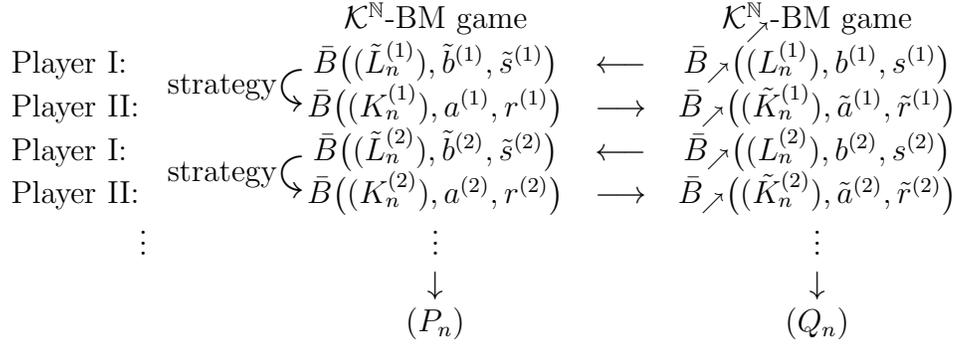

 \centering
 \begin{tabular}{lccc}
  \hspace*{100pt}&$\K^{\N}$-BM game&&$\K_{\nearrow}^{\N}$-BM game\\
  Player~I:&%
  \makebox[0pt][r]{\hspace*{3pt}\raisebox{-8pt}[0pt][0pt]{strategy}\hspace*{2pt}%
                   \raisebox{-15pt}[0pt][0pt]{\rotatebox{90}{\Large$\curvearrowleft$}}\hspace*{4pt}}%
  $\bar{B}\bigl((\tilde{L}_n^{(1)}),\tilde{b}^{(1)},\tilde{s}^{(1)}\bigr)$&$\longleftarrow$&%
  $\bar{B}_{\nearrow}\bigl((L_n^{(1)}),b^{(1)},s^{(1)}\bigr)$\\
  Player~II:&
  $\bar{B}\bigl((K_n^{(1)}),a^{(1)},r^{(1)}\bigr)$&$\longrightarrow$&%
  $\bar{B}_{\nearrow}\bigl((\tilde{K}_n^{(1)}),\tilde{a}^{(1)},\tilde{r}^{(1)}\bigr)$\\
  Player~I:&%
  \makebox[0pt][r]{\hspace*{3pt}\raisebox{-8pt}[0pt][0pt]{strategy}\hspace*{2pt}%
                   \raisebox{-15pt}[0pt][0pt]{\rotatebox{90}{\Large$\curvearrowleft$}}\hspace*{4pt}}%
  $\bar{B}\bigl((\tilde{L}_n^{(2)}),\tilde{b}^{(2)},\tilde{s}^{(2)}\bigr)$&$\longleftarrow$&%
  $\bar{B}_{\nearrow}\bigl((L_n^{(2)}),b^{(2)},s^{(2)}\bigr)$\\
  Player~II:&
  $\bar{B}\bigl((K_n^{(2)}),a^{(2)},r^{(2)}\bigr)$&$\longrightarrow$&%
  $\bar{B}_{\nearrow}\bigl((\tilde{K}_n^{(2)}),\tilde{a}^{(2)},\tilde{r}^{(2)}\bigr)$\\
  \multicolumn{1}{c}{\vdots}&\vdots&&\vdots\\
  &$\downarrow$&&$\downarrow$\\
  &$(P_n)$&&$(Q_n)$
 \end{tabular}
 \caption{Outline of the proof of that $\K^{\N}$-residuality implies $\K_{\nearrow}^{\N}$-residuality}
 \label{fig:transfer_K_nearrow^N-->K^N}
\end{figure}

This subsection is devoted to the outline of the proof
that $\K_{\nearrow}^{\N}$-residuality implies $\K^{\N}$-residuality,
or equivalently, that if the $\K_{\nearrow}^{\N}$-BM game has a winning strategy for Player~II
then so does the $\K^{\N}$-BM game.
Figure~\ref{fig:transfer_K^N-->K_nearrow^N} illustrates this,
and Figure~\ref{fig:transfer_K_nearrow^N-->K^N} allows us
to guess easily the outline of the proof of the other implication.

Suppose that Player~I chose $\bar{B}\bigl((K_n^{(1)}),a^{(1)},r^{(1)}\bigr)$ in the first turn.
Player~II transfers it to a certain set,
say $\bar{B}_{\nearrow}\bigl((\tilde{K}_n^{(1)}),\tilde{a}^{(1)},\tilde{r}^{(1)}\bigr)$,
in the $\K_{\nearrow}^{\N}$-BM game.
Then the winning strategy in the $\K_{\nearrow}^{\N}$-BM game tells Player~II
to take a set $\bar{B}_{\nearrow}\bigl((L_n^{(1)}),b^{(1)},s^{(1)}\bigr)$.
Player~II transfers it to a set $\bar{B}\bigl(\tilde{L}_n^{(1)},\tilde{b}^{(1)},\tilde{s}^{(1)}\bigr)$,
which will be the real reply in the $\K^{\N}$-BM game.
In a similar way, after Player~I replies $\bar{B}\bigl((K_n^{(2)}),a^{(2)},r^{(2)}\bigr)$,
Player~II obtains $\bar{B}_{\nearrow}\bigl((\tilde{K}_n^{(2)}),\tilde{a}^{(2)},\tilde{r}^{(2)}\bigr)$,
$\bar{B}_{\nearrow}\bigl((L_n^{(2)}),b^{(2)},s^{(2)}\bigr)$,
and $\bar{B}\bigl((\tilde{L}_n^{(2)}),\tilde{b}^{(2)},\tilde{s}^{(2)}\bigr)$.
Player~II continues this strategy.

Since $\K^{\N}$ and $\K_{\nearrow}^{\N}$ are compact,
the intersections of the closed sets chosen by the players are nonempty.
By modifying the winning strategy for the $\K_{\nearrow}^{\N}$-BM game,
we may assume that $\lim_{m\to\infty}s^{(m)}=0$,
so that the intersection in this game is a singleton.
Furthermore,
since the transfers are executed so that $\tilde{s}^{(m)}\leqq s^{(m)}$ holds for every $m\in\N$
as will be stated below,
the intersection in the $\K^{\N}$-BM game is also a singleton.

We write
\begin{align*}
 &\bigcap_{m=1}^{\infty}\bar{B}\bigl((K_n^{(m)}),a^{(m)},r^{(m)}\bigr)
 =\bigcap_{m=1}^{\infty}\bar{B}\bigl((\tilde{L}_n^{(m)}),\tilde{b}^{(m)},\tilde{s}^{(m)}\bigr)
  =\bigl\{(P_n)\bigr\}\quad\text{and}\\
 &\bigcap_{m=1}^{\infty}\bar{B}_{\nearrow}\bigl((\tilde{K}_n^{(m)}),\tilde{a}^{(m)},\tilde{r}^{(m)}\bigr)
 =\bigcap_{m=1}^{\infty}\bar{B}_{\nearrow}\bigl((L_n^{(m)}),b^{(m)},s^{(m)}\bigr)
  =\bigl\{(Q_n)\bigr\}.
\end{align*}
Notice that
\[
 \lim_{m\to\infty}(K_n^{(m)})=\lim_{m\to\infty}(\tilde{L}_n^{(m)})=(P_n)\quad\text{and}\quad
 \lim_{m\to\infty}(\tilde{K}_n^{(m)})=\lim_{m\to\infty}(L_n^{(m)})=(Q_n).
\]
Since Player~II follows the winning strategy in the $\K_{\nearrow}^{\N}$-BM game,
we have $(Q_n)\in\K_{\F}^{\N}\cap\K_{\nearrow}^{\N}$,
or equivalently $\bigcup_{n=1}^{\infty}Q_n\in\F$.
Thus all we have to show is that $(P_n)\in\K_{\F}^{\N}$,
and to this aim it suffices to prove that $\bigcup_{n=1}^{\infty}P_n=\bigcup_{n=1}^{\infty}Q_n$.

\subsection{Details of the transfers}
\subsubsection{Conditions and definitions}
A \word{stage} consists of two moves
(one in the $\K^{\N}$-BM game and one in the $\K_{\nearrow}^{\N}$-BM game)
which lie at the same height in Figures~1 and 2.
When we describe the situation at a fixed stage,
we omit the integer $m$ indicating the stage unless ambiguity may be caused:
for example, we write $K_n$ in place of $K_n^{(m)}$.
This is not only for simple notation;
we try to offer explanation of the transfers
which will go in the proofs of both implications,
and this omission solves the problem that
when we describe the stage having, say, $K_n^{(m)}$,
the previous stage can have $L_n^{(m-1)}$ or $L_n^{(m)}$
depending on which implication we look at.

The transfers are executed so that the following conditions,
written as (\textasteriskcentered) afterwards, are fulfilled:
\begin{enumerate}
 \item $\tilde{a}\geqq a$, $\tilde{b}\geqq b$, $\tilde{r}\leqq r/2$, and $\tilde{s}\leqq s/2$;
 \item $\bigcup_{j=1}^{n}K_j\subset\tilde{K}_n$ for $n=1,\ldots,a$, and
       $\bigcup_{j=1}^{n}\tilde{L}_j\subset L_n$ for $n=1,\ldots,b$;\label{enum:affiliation}
 \item $\bigcup_{n=1}^{a}K_n=\tilde{K}_{\tilde{a}}$
       and $\bigcup_{n=1}^{\tilde{b}}\tilde{L}_n=L_b$.
\end{enumerate}

For $x\in\bigcup_{n=1}^{a}K_n=\tilde{K}_{\tilde{a}}$,
its \word{affiliation} $(n_1,n_2)$ is the pair of
the integer $n_1\in\{1,\ldots,a\}$ with $x\in K_{n_1}$, called the \word{first affiliation} of $x$,
and the least integer $n_2\in\{1,\ldots,\tilde{a}\}$ with $x\in\tilde{K}_{n_2}$,
called the \word{second affiliation} of $x$.
We give a similar definition for the points in $\bigcup_{n=1}^{\tilde{b}}\tilde{L}_n=L_b$:
for $x\in\bigcup_{n=1}^{\tilde{b}}\tilde{L}_n=L_b$,
its \word{affiliation} $(n_1,n_2)$ is the pair of
the integer $n_1\in\{1,\ldots,\tilde{b}\}$ with $x\in\tilde{L}_{n_1}$,
called the \word{first affiliation} of $x$,
and the least integer $n_2\in\{1,\ldots,b\}$ with $x\in L_{n_2}$,
called the \word{second affiliation} of $x$.
Strictly speaking, we should specify the stage at which the affiliations are defined,
because, for instance,
it may be that $L_{b^{(m)}}^{(m)}\cap L_{b^{(m')}}^{(m')}\ne\emptyset$ for distinct $m$ and $m'$.
However, since we can easily guess the stage from the context,
we choose not to specify it in order to avoid complexity.

\begin{rem}
 Condition (\ref{enum:affiliation}) in (\textasteriskcentered) is equivalent to the condition that
 the first affiliation is always greater than or equal to the second affiliation.
\end{rem}

Let us look at $\bar{B}\bigl((K_n),a,r\bigr)\in\B$ and
$\bar{B}_{\nearrow}\bigl((\tilde{K}_n),\tilde{a},\tilde{r}\bigr)\in\B_{\nearrow}$
at any stage except the first one.
We have $\bar{B}\bigl((\tilde{L}_n),\tilde{b},\tilde{s}\bigr)\in\B$
and $\bar{B}_{\nearrow}\bigl((L_n),b,s\bigr)\in\B_{\nearrow}$
at the previous stage.
Since $\bar{B}\bigl((K_n),a,r\bigr)\subset\bar{B}\bigl((\tilde{L}_n),\tilde{b},\tilde{s}\bigr)$,
for each $x\in\bigcup_{n=1}^{\tilde{b}}K_n$
there exists a unique $y\in\bigcup_{n=1}^{\tilde{b}}\tilde{L}_n=L_b$
satisfying $\rho(x,y)\leqq\tilde{s}$,
where uniqueness follows from
the assumption that any two distinct points in $\bigcup_{n=1}^{\tilde{b}}\tilde{L}_n$
have distance at least $3\tilde{s}$.
This $y$ is called the \word{parent} of $x$.
Observe that if $x\in K_n$ then $y\in\tilde{L}_n$.
We give a similar definition also when we look at
$\bar{B}_{\nearrow}\bigl((L_n),b,s\bigr)\in\B_{\nearrow}$ and
$\bar{B}\bigl((\tilde{L}_n),\tilde{b},\tilde{s}\bigr)\in\B$:
the \word{parent} of $x\in L_{\tilde{a}}$ is the unique $y\in\bigcup_{n=1}^{a}K_n=\tilde{K}_{\tilde{a}}$
satisfying $\rho(x,y)\leqq\tilde{r}$.

\subsubsection{Transfers from the $\K^{\N}$-BM game to the $\K_{\nearrow}^{\N}$-BM game}
Given a move $\bar{B}\bigl((K_n),a,r\bigr)\in\B$,
we shall construct its transfer
$\bar{B}_{\nearrow}\bigl((\tilde{K}_n),\tilde{a},\tilde{r}\bigr)\in\B_{\nearrow}$.
If it is the first move of Player I,
then we put $\tilde{a}=a$, $\tilde{r}=r/2$, and
$\tilde{K}_n=\bigcup_{j=1}^{n}K_j$ for every $n\in\N$,
and we can easily see that the conditions (\textasteriskcentered) are fulfilled.
So suppose otherwise.
Then we already know $\bar{B}_{\nearrow}\bigl((L_n),b,s\bigr)\in\B_{\nearrow}$
and its transfer $\bar{B}\bigl((\tilde{L}_n),\tilde{b},\tilde{s}\bigr)\in\B$,
and we have $\bar{B}\bigl((K_n),a,r\bigr)\subset\bar{B}\bigl((\tilde{L}_n),\tilde{b},\tilde{s}\bigr)$.

Put $\tilde{a}=a$ and $\tilde{r}=\min\{s-\tilde{s},r/2\}$, and
define $\tilde{K}_n=\bigcup_{j=1}^{n}K_j$ for $n>\tilde{b}$.
We define $\tilde{K}_n$ for $n\leqq\tilde{b}$
by declaring that the second affiliation of each $x\in\bigcup_{n=1}^{\tilde{b}}K_n$
is the same as that of the parent of $x$.

\begin{claim}
 We have $d(\tilde{K}_n,L_n)\leqq\tilde{s}$ for $n=1,\ldots,b$.
\end{claim}

\begin{prf}
 Fix such an integer $n$.

 Let $x\in\tilde{K}_n$ and denote its affiliation by $(n_1,n_2)$.
 Then the parent $y$ of $x$ has affiliation $(n_1,n_2)$ and so belongs to $L_{n_2}$.
 It follows from $y\in L_{n_2}\subset L_n$ and $\rho(x,y)\leqq\tilde{s}$
 that $x\in L_n[\tilde{s}]$.

 Conversely let $y\in L_n$ and denote its affiliation by $(n_1,n_2)$.
 Then there exists a point $x\in K_{n_1}$ with $\rho(x,y)\leqq\tilde{s}$
 because $d(K_{n_1},\tilde{L}_{n_1})\leqq\tilde{s}$.
 Since $y$ is the parent of $x$,
 the affiliation of $x$ is $(n_1,n_2)$.
 Therefore $x\in\tilde{K}_{n_2}\subset\tilde{K}_n$ and so $y\in\tilde{K}_n[\tilde{s}]$.
\end{prf}

We may deduce from this claim that
$\bar{B}_{\nearrow}\bigl((\tilde{K}_n),\tilde{a},\tilde{r}\bigr)
 \subset\bar{B}_{\nearrow}\bigl((L_n),b,s\bigr)$
using the triangle inequality and $\tilde{r}+\tilde{s}\leqq s$.
Therefore $\bar{B}_{\nearrow}\bigl((\tilde{K}_n),\tilde{a},\tilde{r}\bigr)$
is a valid reply in the $\K_{\nearrow}^{\N}$-BM game.
It is easy to see that the conditions (\textasteriskcentered) are fulfilled.

\subsubsection{Transfers from the $\K_{\nearrow}^{\N}$-BM game to the $\K^{\N}$-BM game}
Given a move $\bar{B}_{\nearrow}\bigl((L_n),b,s\bigr)\in\B_{\nearrow}$,
we shall construct its transfer $\bar{B}\bigl((\tilde{L}_n),\tilde{b},\tilde{s}\bigr)\in\B$.
If it is the first move of Player~I,
then we put $\tilde{b}=b$, $\tilde{s}=s/2$,
$\tilde{L}_1=L_1$, and $\tilde{L}_n=L_n\setminus L_{n-1}$ for every $n\geqq2$.
We can easily see that the conditions (\textasteriskcentered) are fulfilled in this case.
So suppose otherwise.
Then we already know $\bar{B}\bigl((K_n),a,r\bigr)\in\B$
and its transfer $\bar{B}_{\nearrow}\bigl((\tilde{K}_n),\tilde{a},\tilde{r}\bigr)\in\B_{\nearrow}$,
and we have
$\bar{B}_{\nearrow}\bigl((L_n),b,s\bigr)
 \subset\bar{B}_{\nearrow}\bigl((\tilde{K}_n),\tilde{a},\tilde{r}\bigr)$.

Put $\tilde{b}=b+1$ and $\tilde{s}=\min\{r-\tilde{r},s/2\}$,
and define $\tilde{L}_n=L_{n-1}$ for $n>\tilde{b}$.
We define $\tilde{L}_n$ for $n\leqq\tilde{b}$
by determining the first affiliation of each point in $L_b$ as follows.
Let $x\in L_b$ and denote its second affiliation by $n_2$.
If $n_2>\tilde{a}$, then the first affiliation of $x$ is $n_2$.
Suppose $n_2\leqq\tilde{a}$, and let $y\in\tilde{K}_{n_2}$ denote the parent of $x$.
If the second affiliation of $y$ is $n_2$,
then the first affiliation of $x$ is the same as that of $y$;
otherwise the first affiliation of $x$ is $\tilde{b}$.

\begin{claim}
 We have $d(\tilde{L}_n,K_n)\leqq\tilde{r}$ for $n=1,\ldots,a$.
\end{claim}

\begin{prf}
 Fix such an integer $n$.

 Let $x\in\tilde{L}_n$ and denote its parent by $y$.
 Then it follows that $x$ and $y$ have the same affiliation, and so $y\in K_n$.
 Hence we may infer from $\rho(x,y)\leqq\tilde{r}$ that $x\in K_n[\tilde{r}]$.

 Conversely let $y\in K_n$ and denote its second affiliation by $n_2$.
 Then there exists a point $x\in L_{n_2}$ with $\rho(x,y)\leqq\tilde{r}$
 because $d(\tilde{K}_{n_2},L_{n_2})\leqq\tilde{r}$.
 Since $y$ is the parent of $x$ and has the same second affiliation as $x$,
 the first affiliation of $x$ is $n$.
 Therefore $y\in\tilde{L}_n[\tilde{r}]$.
\end{prf}

We may deduce from the claim that
$\bar{B}\bigl((\tilde{L}_n),\tilde{b},\tilde{s}\bigr)\subset\bar{B}\bigl((K_n),a,r\bigr)$
using the triangle inequality and $\tilde{r}+\tilde{s}\leqq r$.
Therefore $\bar{B}\bigl((\tilde{L}_n),\tilde{b},\tilde{s}\bigr)$ is a valid reply
in the $\K^{\N}$-BM game.
It is easy to see that the conditions (\textasteriskcentered) are fulfilled.

\subsection{Proof of $\bigcup_{n=1}^{\infty}P_n=\bigcup_{n=1}^{\infty}Q_n$}
We shall prove that $\bigcup_{n=1}^{\infty}P_n=\bigcup_{n=1}^{\infty}Q_n$,
which will complete the proof of Theorem~\ref{thm:main-BM} and hence of our main theorem.
Recall that $(K_n^{(m)})$ and $(\tilde{K}_n^{(m)})$ converge to $(P_n)$ and $(Q_n)$
respectively as $m$ tends to infinity.
In other words we have $\lim_{m\to\infty}K_n^{(m)}=P_n$ and $\lim_{m\to\infty}\tilde{K}_n^{(m)}=Q_n$
for every $n\in\N$.

In order to prove $\bigcup_{n=1}^{\infty}P_n\subset\bigcup_{n=1}^{\infty}Q_n$,
it is enough to show that $\bigcup_{j=1}^{n}P_j\subset Q_n$ for every $n\in\N$.
The set $\bigl\{\,(A,B)\in\K^2\bigm|A\subset B\,\bigr\}$ is closed in $\K^2$
and contains $(\bigcup_{j=1}^{n}K_j^{(m)},\tilde{K}_n^{(m)})$ for all $m\in\N$.
Since $(\bigcup_{j=1}^{n}K_j^{(m)},\tilde{K}_n^{(m)})$ converges to
$(\bigcup_{j=1}^{n}P_j,Q_n)$ as $m$ tends to infinity,
which follows from the continuity of the map
$(A_1,\ldots,A_n)\longmapsto\bigcup_{j=1}^{n}A_j$ from $\K^n$ to $\K$,
we obtain $\bigcup_{j=1}^{n}P_j\subset Q_n$.

Now we shall prove $\bigcup_{n=1}^{\infty}Q_n\subset\bigcup_{n=1}^{\infty}P_n$.
Let $x\in\bigcup_{n=1}^{\infty}Q_n$,
and denote by $n$ the least positive integer with $x\in Q_n$.
Since it is easy to observe that $K_1^{(m)}=\tilde{K}_1^{(m)}$ for every $m\in\N$,
which implies $P_1=Q_1$,
we may assume that $n\geqq2$.
Because $Q_{n-1}$ is closed and $x\notin Q_{n-1}$,
there exists a positive real number $r$ less than $1$
satisfying $\bar{B}(x,4r)\cap Q_{n-1}=\emptyset$, that is, $x\notin Q_{n-1}[4r]$.
Fix a positive integer $m_0$ such that
$\tilde{a}^{(m)}\geqq n$, $\tilde{r}^{(m)}\leqq r$,
and $d(\tilde{K}_{n-1}^{(m)},Q_{n-1})\leqq r$
for every $m\geqq m_0$.
Observe that $x\notin\tilde{K}_{n-1}^{(m)}[3r]$ for every $m\geqq m_0$.

Set $k_0=\lceil1/r\rceil$.
For each $k\geqq k_0$,
choose $m_k\geqq m_0$ satisfying $d(\tilde{K}_n^{(m)},Q_n)\leqq1/k$ for every $m\geqq m_k$,
and for each $m\geqq m_k$ take $y_{km}\in\tilde{K}_n^{(m)}$
with $\rho(x,y_{km})\leqq1/k$
and let $z_{km}\in\tilde{K}_{n}^{(m_0)}$ denote the unique point
satisfying $\rho(y_{km},z_{km})\leqq\tilde{r}^{(m_0)}$.

\begin{claim}
 The two points $y_{km}$ and $z_{km}$ have the same affiliation.
\end{claim}

\begin{prf}
 By an \word{ancestor} of $y_{km}$ we mean a point that can be written
 as `the parent of \ldots\ the parent of $y_{km}$.'
 Observe that $z_{km}$ is an ancestor of $y_{km}$.
 Indeed if we denote by $z_{km}'$ the ancestor of $y_{km}$ in $\tilde{K}_n^{(m_0)}$, then
 \[
  \rho(y_{km},z_{km}')<\tilde{r}^{(m_0)}+\frac{\tilde{r}^{(m_0)}}{2}+\frac{\tilde{r}^{(m_0)}}{2^2}+\cdots=2\tilde{r}^{(m_0)}
 \]
 and so $\rho(z_{km},z_{km}')<3\tilde{r}^{(m_0)}$,
 which implies $z_{km}=z_{km}'$.

 In order to prove our claim,
 it suffices to prove that the second affiliation of
 the ancestor $w\in\tilde{K}_n^{(m')}$ of $y_{km}$ is $n$ for any $m'\in\{m_0,\ldots,m\}$.
 We can see $\rho(w,y_{km})\leqq2\tilde{r}^{(m')}\leqq2r$ by the same reasoning as above.
 Therefore we have
 \[
  \rho(w,x)\leqq\rho(w,y_{km})+\rho(y_{km},x)\leqq2r+\frac{1}{k}\leqq3r.
 \]
 Thus the second affiliation of $w$ cannot be less than $n$ because $x\notin\tilde{K}_{n-1}^{(m')}[3r]$.
\end{prf}

Note that all $z_{km}$ belong to the single finite set $\tilde{K}_n^{(m_0)}$.
We can choose $z_k\in K_n^{(m_0)}$ for $k\geqq k_0$ inductively so that
the set
\[
 \{\,m\geqq m_k\mid z_{k_0m}=z_{k_0},\ldots,z_{km}=z_k\,\}
\]
is infinite for any $k\geqq k_0$.
Then we take $z\in K_n^{(m_0)}$ for which $\{\,k\geqq k_0\mid z_k=z\,\}$ is infinite,
and put $\{\,k\geqq k_0\mid z_k=z\,\}=\{k_1,k_2,\ldots\}$, where $k_1<k_2<\cdots$.
Since the set
\[
 \{\,m\geqq m_{k_j}\mid z_{k_1m}=\cdots=z_{k_jm}=z\,\}
\]
is infinite for every $j\in\N$,
we may construct a strictly increasing sequence $m_1'$, $m_2'$, \ldots\ of positive integers
satisfying $z_{k_1m_j'}=\cdots=z_{k_jm_j'}=z$.

Let $l$ denote the first affiliation of $z$.
Then the foregoing claim shows that whenever $i\leqq j$,
the first affiliation of $y_{k_im_j'}$ is $l$,
which implies that $x\in K_l^{(m_j')}[1/k_i]$.
For any $i\in\N$,
since $d(K_l^{(m_j')},P_l)\leqq1/k_i$ for sufficiently large $j$,
we have $x\in P_l[2/k_i]$.
Hence $x\in\bigcap_{i=1}^{\infty}P_l[2/k_i]=P_l$.
This completes the proof.

\end{document}